# Normal approximation for nonlinear statistics using a concentration inequality approach

LOUIS H.Y. CHEN[1] and QI-MAN SHAO[2]

[1] *Institute for Mathematical Sciences, National University of Singapore, Singapore 118402, Republic of Singapore. E-mail: matchyl@nus.edu.sg*
[2] *Department of Mathematics, Hong Kong University of Science and Technology, Clear Water Bay, Kowloon, Hong Kong, China; Department of Mathematics, University of Oregon, Eugene, OR 97403, USA; Department of Mathematics, Zhejiang University, Hangzhou, Zhejiang 310027, China. E-mail: maqmshao@ust.hk*

Let $T$ be a general sampling statistic that can be written as a linear statistic plus an error term. Uniform and non-uniform Berry–Esseen type bounds for $T$ are obtained. The bounds are the best possible for many known statistics. Applications to U-statistics, multisample U-statistics, L-statistics, random sums and functions of nonlinear statistics are discussed.

## 1. Introduction

Let $X_1, X_2, \ldots, X_n$ be independent random variables and let $T := T(X_1, \ldots, X_n)$ be a general sampling statistic. In many cases $T$ can be written as a linear statistic plus an error term, say $T = W + \Delta$, where

$$W = \sum_{i=1}^n g_i(X_i), \qquad \Delta := \Delta(X_1, \ldots, X_n) = T - W$$

and $g_i := g_{n,i}$ are Borel measurable functions. Typical cases include U-statistics, multisample U-statistics, L-statistics and random sums. Assume that

$$\mathrm{E}(g_i(X_i)) = 0 \qquad \text{for } i = 1, 2, \ldots, n \quad \text{and} \quad \sum_{i=1}^n \mathrm{E}(g_i^2(X_i)) = 1. \tag{1.1}$$

It is clear that if $\Delta \to 0$ in probability as $n \to \infty$, then we have the central limit theorem

$$\sup_z |P(T \le z) - \Phi(z)| \to 0, \tag{1.2}$$

where $\Phi$ denotes the standard normal distribution function, provided that the Lindeberg condition holds:

$$\forall \varepsilon > 0, \qquad \sum_{i=1}^n \mathrm{E} g_i^2(X_i) I(|g_i(X_i)| > \varepsilon) \to 0.$$

If in addition, $\mathrm{E}|\Delta|^p < \infty$ for some $p > 0$, then by the Chebyshev inequality, one can obtain the rate of convergence

$$\sup_z |P(T \le z) - \Phi(z)| \le \sup_z |P(W \le z) - \Phi(z)| + 2(\mathrm{E}|\Delta|^p)^{1/(1+p)}. \tag{1.3}$$

The first term on the right-hand side of (1.3) is well understood via the Berry–Esseen inequality. For example, using Stein's method, Chen and Shao [9] obtained

$$\begin{aligned} &\sup_z |P(W \le z) - \Phi(z)| \\ &\quad \le 4.1 \left( \sum_{i=1}^n \mathrm{E} g_i^2(X_i) I(|g_i(X_i)| > 1) + \sum_{i=1}^n \mathrm{E}|g_i(X_i)|^3 I(|g_i(X_i)| \le 1) \right). \end{aligned} \tag{1.4}$$

However, the bound $(\mathrm{E}|\Delta|^p)^{1/(1+p)}$ is, in general, not sharp for many commonly used statistics. Many authors have worked toward obtaining better Berry–Esseen bounds. For example, sharp Berry–Esseen bounds have been obtained for general symmetric statistics by van Zwet [24] and Friedrich [12]. An Edgeworth expansion with remainder $O(n^{-1})$ for symmetric statistics was proved by Bentkus, Götze and van Zwet [3].

The main purpose of this paper is to establish uniform and non-uniform Berry–Esseen bounds for general nonlinear statistics. The bounds are the best possible for many known statistics. Our proof is based on a randomized concentration inequality approach to bounding $P(W + \Delta \le z) - P(W \le z)$. Because proofs of uniform and non-uniform bounds for sums of independent random variables can be proved via Stein's method (Chen and Shao [9]), which is much neater and simpler than the traditional Fourier analysis approach, this paper provides a direct and unifying treatment toward the Berry–Esseen bounds for general nonlinear statistics.

This paper is organized as follows. The main results are stated in next section, three applications are presented in Section 3 and an example is given in Section 4 to show the sharpness of the main results. Proofs of the main results are given in Section 5. For the proofs of other results, including Example 4.1, the reader is referred to our technical report (Chen and Shao [10]).

Throughout this paper, $C$ will denote an absolute constant whose value may change at each appearance. The $L_p$ norm of a random variable $X$ is denoted by $\|X\|_p$, that is, $\|X\|_p = (\mathrm{E}|X|^p)^{1/p}$ for $p \ge 1$.

## 2. Main results

Let $\{X_i, 1 \le i \le n\}$, $T$, $W$ and $\Delta$ be defined as in Section 1. In the following theorems, we assume that (1.1) is satisfied. Put

$$\beta = \sum_{i=1}^n \mathrm{E}|g_i(X_i)|^2 I(|g_i(X_i)| > 1) + \sum_{i=1}^n \mathrm{E}|g_i(X_i)|^3 I(|g_i(X_i)| \le 1) \tag{2.1}$$

and let $\delta > 0$ satisfy

$$\sum_{i=1}^n \mathrm{E}|g_i(X_i)| \min(\delta, |g_i(X_i)|) \ge 1/2. \tag{2.2}$$

**Theorem 2.1.** *For each $1 \le i \le n$, let $\Delta_i$ be a random variable such that $X_i$ and $(\Delta_i, W - g_i(X_i))$ are independent. Then*

$$\sup_z |P(T \le z) - P(W \le z)| \le 4\delta + \mathrm{E}|W\Delta| + \sum_{i=1}^n \mathrm{E}|g_i(X_i)(\Delta - \Delta_i)| \tag{2.3}$$

*for $\delta$ satisfying (2.2). In particular, we have*

$$\sup_z |P(T \le z) - P(W \le z)| \le 2\beta + \mathrm{E}|W\Delta| + \sum_{i=1}^n \mathrm{E}|g_i(X_i)(\Delta - \Delta_i)| \tag{2.4}$$

*and*

$$\sup_z |P(T \le z) - \Phi(z)| \le 6.1\beta + \mathrm{E}|W\Delta| + \sum_{i=1}^n \mathrm{E}|g_i(X_i)(\Delta - \Delta_i)|. \tag{2.5}$$

The next theorem provides a non-uniform bound.

**Theorem 2.2.** *For each $1 \le i \le n$, let $\Delta_i$ be a random variable such that $X_i$ and $(\Delta_i, \{X_j, j \ne i\})$ are independent. Then, for $\delta$ satisfying (2.2) and for $z \in R^1$,*

$$|P(T \le z) - P(W \le z)| \le \gamma_z + \mathrm{e}^{-|z|/3}\tau, \tag{2.6}$$

*where*

$$\gamma_z = P(|\Delta| > (|z|+1)/3) + \sum_{i=1}^n P(|g_i(X_i)| > (|z|+1)/3) + \sum_{i=1}^n P(|W - g_i(X_i)| > (|z|-2)/3) P(|g_i(X_i)| > 1), \tag{2.7}$$

$$\tau = 22\delta + 8.5\|\Delta\|_2 + 3.6\sum_{i=1}^n \|g_i(X_i)\|_2\|\Delta - \Delta_i\|_2. \tag{2.8}$$

*In particular, if* $\mathrm{E}|g_i(X_i)|^p < \infty$ *for* $2 < p \le 3$*, then*

$$\begin{aligned}
&|P(T \le z) - \Phi(z)| \\
&\quad \le P(|\Delta| > (|z|+1)/3) \\
&\qquad + C(|z|+1)^{-p}\left(\|\Delta\|_2 + \sum_{i=1}^n \|g_i(X_i)\|_2\|\Delta - \Delta_i\|_2 + \sum_{i=1}^n \mathrm{E}|g_i(X_i)|^p\right).
\end{aligned} \tag{2.9}$$

A result similar to (2.5) was obtained by Friedrich [12] for $g_i = \mathrm{E}(T|X_i)$ using the method of characteristic function. Our proof is direct and simpler, and the bounds are easier to calculate. The non-uniform bounds in (2.6) and (2.9) for general nonlinear statistics are new.

***Remark 2.1.*** Assume $\mathrm{E}|g_i(X_i)|^p < \infty$ for $p > 2$. Let

$$\delta = \left(\frac{2(p-2)^{p-2}}{(p-1)^{p-1}}\sum_{i=1}^n \mathrm{E}|g_i(X_i)|^p\right)^{1/(p-2)}. \tag{2.10}$$

Then (2.2) is satisfied. This follows from the inequality

$$\min(a,b) \ge a - \frac{(p-2)^{p-2}a^{p-1}}{(p-1)^{p-1}b^{p-2}} \tag{2.11}$$

for $a \ge 0$ and $b > 0$.

***Remark 2.2.*** If $\beta \le 1/2$, then (2.2) is satisfied with $\delta = \beta/2$.

***Remark 2.3.*** Let $\delta > 0$ be such that

$$\sum_{i=1}^n \mathrm{E}g_i^2(X_i)I(|g_i(X_i)| > \delta) \le 1/2.$$

Then (2.2) holds. In particular, if $X_1, X_2, \ldots, X_n$ are independent and identically distributed (i.i.d.) random variables and $g_i = g_1$, then (2.2) is satisfied with $\delta = c_0/\sqrt{n}$, where $c_0$ is a constant such that $\mathrm{E}(\sqrt{n}g_1(X_1))^2 I(|\sqrt{n}g_1(X_1)| > c_0) \le 1/2$.

***Remark 2.4.*** In Theorems 2.1 and 2.2, the choice of $\Delta_i$ is flexible. For example, one can choose $\Delta_i = \Delta(X_1, \ldots, X_{i-1}, 0, X_{i+1}, \ldots, X_n)$ or $\Delta_i = \Delta(X_1, \ldots, X_{i-1}, \hat{X}_i, X_{i+1}, \ldots, X_n)$, where $\{\hat{X}_i, 1 \le i \le n\}$ is an independent copy of $\{X_i, 1 \le i \le n\}$. The choice of $g_i$ is also flexible. It can be more general than $g_i(x) = \mathrm{E}(T|X_i = x)$, which is commonly used in the literature.

***Remark 2.5.*** Let $X_1, \ldots, X_n$ be independent normally distributed random variables with mean zero and variance $1/n$, and let $W$, $T$ and $\Delta$ be as in Example 4.1. Then

$$\begin{aligned} & \mathrm{E}|W\Delta| + \sum_{i=1}^{n} \mathrm{E}|X_i|^3 + \sum_{i=1}^{n} \mathrm{E}|X_i(\Delta(X_1, \ldots, X_i, \ldots, X_n) - \Delta(X_1, \ldots, 0, \ldots, X_n))| \\ & \leq C\varepsilon^{2/3} \end{aligned} \tag{2.12}$$

for $(1/\varepsilon)^{4/3} \leq n \leq 16(1/\varepsilon)^{4/3}$. This together with (4.5) shows that the bound in (2.4) is achievable. Moreover, the term $\sum \mathrm{E}|g_i(X_i)(\Delta - \Delta_i)|$ in (2.4) cannot be dropped.

# 3. Applications

Theorems 2.1 and 2.2 can be applied to a wide range of different statistics and provide the bounds of the best possible order in many instances. To illustrate the usefulness and the generality of these results, we give three applications in this section. The uniform bounds refine many existing results by specifying absolute constants, while the non-uniform bounds are new for many cases. Because the results are direct applications of Theorems 2.1 and 2.2, and the proofs are more or less routine verifications of assumptions, the proofs are omitted. Refer to Chen and Shao [10] for detailed proofs. One may also refer to Chen and Shao [10] for applications to random sums of independent random variables with non-random centering and to functions of nonlinear statistics.

## 3.1. U-statistics

Let $X_1, X_2, \ldots, X_n$ be a sequence of i.i.d. random variables, and let $h(x_1, \ldots, x_m)$ be a real-valued Borel measurable symmetric function of $m$ variables, where $m$ $(2 \leq m < n)$ may depend on $n$. Consider the Hoeffding [17] U-statistic

$$U_n = \binom{n}{m}^{-1} \sum_{1 \leq i_1 < \cdots < i_m \leq n} h(X_{i_1}, \ldots, X_{i_m}).$$

The U-statistic elegantly and usefully generalizes the notion of a sample mean. Numerous investigations on the limiting properties of the U-statistics have been done during the last few decades. A systematic presentation of the theory of U-statistics was given by Koroljuk and Borovskich [19]. We refer to the studies on uniform Berry–Esseen bounds for U-statistics by Filippova [11], Grams and Serfling [13], Bickel [4], Chan and Wierman [8], Callaert and Janssen [7], Serfling [21], van Zwet [24] and Friedrich [12]. One can also refer to the work of Wang, Jing and Zhao [26] on uniform Berry–Esseen bounds for studentized U-statistics.

Applying Theorems 2.1 and 2.2 to the U-statistic yields the following result:

**Theorem 3.1.** *Assume that* $\mathrm{E}h(X_1,\dots,X_m)=0$ *and* $\sigma^2=\mathrm{E}h^2(X_1,\dots,X_m)<\infty$. *Let* $g(x)=\mathrm{E}(h(X_1,X_2,\dots,X_m)|X_1=x)$ *and* $\sigma_1^2=\mathrm{E}g^2(X_1)$. *Assume that* $\sigma_1>0$. *Then*

$$\sup_z\left|P\left(\frac{\sqrt{n}}{m\sigma_1}U_n\le z\right)-P\left(\frac{1}{\sqrt{n}\sigma_1}\sum_{i=1}^n g(X_i)\le z\right)\right|\le\frac{(1+\sqrt{2})(m-1)\sigma}{(m(n-m+1))^{1/2}\sigma_1}+\frac{c_0}{\sqrt{n}},\quad(3.1)$$

*where* $c_0$ *is a constant such that* $\mathrm{E}g^2(X_1)I(|g(X_1)|>c_0\sigma_1)\le\sigma_1^2/2$. *If in addition* $\mathrm{E}|g(X_1)|^p<\infty$ *for* $2<p\le 3$, *then*

$$\sup_z\left|P\left(\frac{\sqrt{n}}{m\sigma_1}U_n\le z\right)-\Phi(z)\right|\le\frac{(1+\sqrt{2})(m-1)\sigma}{(m(n-m+1))^{1/2}\sigma_1}+\frac{6.1\mathrm{E}|g(X_1)|^p}{n^{(p-2)/2}\sigma_1^p}\qquad(3.2)$$

*and, for* $z\in R^1$,

$$\left|P\left(\frac{\sqrt{n}}{m\sigma_1}U_n\le z\right)-\Phi(z)\right|\le\frac{9m\sigma^2}{(1+|z|)^2(n-m+1)\sigma_1^2}+\frac{13.5\mathrm{e}^{-|z|/3}m^{1/2}\sigma}{(n-m+1)^{1/2}\sigma_1}+\frac{C\mathrm{E}|g(X_1)|^p}{(1+|z|)^p n^{(p-2)/2}\sigma_1^p}.\qquad(3.3)$$

*Moreover, if* $\mathrm{E}|h(X_1,\dots,X_m)|^p<\infty$ *for* $2<p\le 3$, *then for* $z\in R^1$,

$$\left|P\left(\frac{\sqrt{n}}{m\sigma_1}U_n\le z\right)-\Phi(z)\right|\le\frac{Cm^{1/2}\mathrm{E}|h(X_1,\dots,X_m)|^p}{(1+|z|)^p(n-m+1)^{1/2}\sigma_1^p}+\frac{C\mathrm{E}|g(X_1)|^p}{(1+|z|)^p n^{(p-2)/2}\sigma_1^p}.\qquad(3.4)$$

Note that the error in (3.1) is $O(n^{-1/2})$ only under the assumption of the finite second moment of $h$. The result appears to be previously unknown. The uniform bound given in (3.2) is not new, but the specifying constant for general $m$ is new. The finite second moment of $h$ is not the weakest assumption for the uniform bound. Friedrich [12] obtained an $O(n^{-1/2})$ when $\mathrm{E}|h|^{5/3}<\infty$. We refer the reader to Bentkus, Götze and Zitikis [2] and Jing and Zhou [18] for discussions on the necessity of the moment condition.

For the non-uniform bound, Zhao and Chen [27] proved that if $m=2$, $\mathrm{E}|h(X_1,X_2)|^3<\infty$, then

$$\left|P\left(\frac{\sqrt{n}}{m\sigma_1}U_n\le z\right)-\Phi(z)\right|\le An^{-1/2}(1+|z|)^{-3}\qquad(3.5)$$

for $z\in R^1$, where the constant $A$ depends not on $n$ and $z$, but on the moment of $h$. Clearly, (3.4) refines Zhao and Chen's result specifying the relationship of the constant $A$ with the moment condition. After we finished proving Theorem 3.1, Wang [25] informed the second author that he also obtained (3.4) for $m=2$ and $p=3$.

Theorems 2.1 and 2.2 can also be applied to derive uniform and non-uniform bounds for the non-i.i.d. case. Especially the uniform bound of Alberink [1] can be easily recovered.

***Remark 3.1.*** Equation (3.3) implies that

$$\left|P\left(\frac{\sqrt{n}}{m\sigma_1}U_n \le z\right) - \Phi(z)\right| \le \frac{Cm^{1/2}\sigma^2}{(1+|z|)^3(n-m+1)^{1/2}\sigma_1^2} + \frac{C\mathrm{E}|g(X_1)|^p}{(1+|z|)^p n^{(p-2)/2}\sigma_1^p} \tag{3.6}$$

for $|z| \le ((n-m+1)/m)^{1/2}$. For $|z| > ((n-m+1)/m)^{1/2}$, a bound like (3.6) can be easily obtained by using the Chebyshev inequality. On the other hand, if (3.6) holds for any $z \in R^1$, then it appears necessary to assume $\mathrm{E}|h(X_1,\ldots,X_m)|^p < \infty$.

***Remark 3.2.*** Theorem 3.1 shows that the central limit theorem can still hold even when $m = m(n) \to \infty$. However, the explicit Berry–Esseen bounds in (3.1) and (3.2) may not be optimal in terms of $m$. For the special case of elementary symmetric polynomials (i.e., $h(x_1,x_2,\ldots,x_m) = x_1 x_2 \ldots x_m - \mu^m$, where $\mu = \mathrm{E}(X_1) \ne 0$), discussed in van Es and Helmers [23], $\sigma/\sigma_1$ in (3.1) and (3.2) is of order $r^m$, where $r > 1$, which is equivalent to requiring $m \le c\log n$ for some $c > 0$ so that the bound tends to 0, whereas in the result of van Es and Helmers [23], the bound tends to 0 as long as $m = o(n^{1/2})$. The Berry–Esseen bound obtained by van Es and Helmers [23] gives an optimal dependence on $m$.

## 3.2. Multisample U-statistics

Consider $k$ independent sequences $\{X_{j1},\ldots,X_{jn_j}\}$ of i.i.d. random variables, $j = 1,\ldots,k$. Let $h(x_{jl}, l = 1,\ldots,m_j, j = 1,\ldots,k)$ be a measurable function symmetric with respect to $m_j$ arguments of the $j$th set, $m_j \ge 1$, $j = 1,\ldots,k$. Let

$$\theta = \mathrm{E}h(X_{jl}, l = 1,\ldots,m_j, j = 1,\ldots,k).$$

The multisample U-statistic is defined as

$$U_{\bar{n}} = \left\{\prod_{j=1}^{k}\binom{n_j}{m_j}^{-1}\right\}\sum h(X_{jl}, l = i_{j1},\ldots,i_{jm_j}, j = 1,\ldots,k),$$

where $\bar{n} = (n_1,\ldots,n_k)$ and the summation is carried out over all $1 \le i_{j1} < \cdots < i_{jm_j} \le n_j$, $n_j \ge 2m_j$, $j = 1,\ldots,k$. Clearly, $U_{\bar{n}}$ is an unbiased estimate of $\theta$. The two-sample Wilcoxon statistic and the two-sample $\omega^2$-statistic are two typical examples of multisample U-statistics. Without loss of generality, assume $\theta = 0$. For $j = 1,\ldots,k$, define

$$h_j(x) = \mathrm{E}(h(X_{11},\ldots,X_{1m_1};\ldots;X_{k1},\ldots,X_{km_k})|X_{j1} = x)$$

and let $\sigma_j^2 = \mathrm{E}h_j^2(X_{j1})$ and

$$\sigma_{\bar{n}}^2 = \sum_{j=1}^{k}\frac{m_j^2}{n_j}\sigma_j^2.$$

A uniform Berry–Esseen bound with $O((\min_{1\le j\le k} n_j)^{-1/2})$ for the multisample U-statistics was obtained by Helmers and Janssen [15] and Borovskich [6] (see Koroljuk and Borovskich [19], pages 304–311). The next theorem refines their results.

**Theorem 3.2.** *Assume that* $\theta = 0$, $\sigma^2 := \mathrm{E}h^2(X_{11},\dots,X_{1m_1};\dots;X_{k1},\dots,X_{km_k}) < \infty$ *and* $\max_{1\le j\le k}\sigma_j > 0$. *Then for* $2 < p \le 3$,

$$\sup_z |P(\sigma_{\bar n}^{-1}U_{\bar n} \le z) - \Phi(z)| \le \frac{(1+\sqrt{2})\sigma}{\sigma_{\bar n}}\sum_{j=1}^k \frac{m_j^2}{n_j} + \frac{6.6}{\sigma_{\bar n}^p}\sum_{j=1}^k \frac{m_j^p}{n_j^{p-1}}\mathrm{E}|h_j(X_{j1})|^p \tag{3.7}$$

*and for* $z \in R^1$,

$$|P(\sigma_{\bar n}^{-1}U_n \le z) - \Phi(z)| \le \frac{9\sigma^2}{(1+|z|)^2\sigma_{\bar n}^2}\left(\sum_{j=1}^k \frac{m_j^2}{n_j}\right)^2 + 13.5\mathrm{e}^{-|z|/3}\frac{\sigma}{\sigma_{\bar n}}\sum_{j=1}^k \frac{m_j^2}{n_j} + \frac{C}{(1+|z|)^p\sigma_{\bar n}^p}\sum_{j=1}^k \frac{m_j^p \mathrm{E}|h_j(X_{j1})|^p}{n_j^{p-1}}. \tag{3.8}$$

## 3.3. L-statistics

Let $X_1,\dots,X_n$ be i.i.d. random variables with a common distribution function $F$ and let $F_n$ be the empirical distribution function defined by

$$F_n(x) = n^{-1}\sum_{i=1}^n I(X_i \le x), \qquad \text{for } x \in R^1.$$

Let $J(t)$ be a real-valued function on $[0,1]$ and define

$$T(G) = \int_{-\infty}^{\infty} xJ(G(x))\,\mathrm{d}G(x)$$

for non-decreasing measurable function $G$. Put

$$\sigma^2 = \int_{-\infty}^{\infty}\int_{-\infty}^{\infty} J(F(s))J(F(t))F(\min(s,t))(1-F(\max(s,t)))\,\mathrm{d}s\,\mathrm{d}t$$

and

$$g(x) = \int_{-\infty}^{\infty} (I(x\le s) - F(s))J(F(s))\,\mathrm{d}s.$$

The statistic $T(F_n)$ is called an L-statistic (see Serfling [21], Chapter 8). Uniform Berry–Esseen bounds for the L-statistic that smooths $J$ were given by Helmers [14] and Helmers, Janssen and Serfling [16]. Applying Theorems 2.1 and 2.2 yields the following uniform and non-uniform bounds for L-statistics.

**Theorem 3.3.** *Let $n \geq 4$. Assume that $\mathrm{E}X_1^2 < \infty$ and $\mathrm{E}|g(X_1)|^p < \infty$ for $2 < p \leq 3$. If the weight function $J(t)$ is Lipschitz of order 1 on $[0,1]$, that is, there exists a constant $c_0$ such that*

$$|J(t) - J(s)| \leq c_0|t-s|, \qquad \text{for } 0 \leq s,t \leq 1, \tag{3.9}$$

*then*

$$\sup_z |P(\sqrt{n}\sigma^{-1}(T(F_n) - T(F)) \leq z) - \Phi(z)| \leq \frac{(1+\sqrt{2})c_0\|X_1\|_2}{\sqrt{n}\sigma} + \frac{6.1\mathrm{E}|g(X_1)|^p}{n^{(p-2)/2}\sigma^p} \tag{3.10}$$

*and*

$$\begin{aligned} &|P(\sqrt{n}\sigma^{-1}(T(F_n) - T(F)) \leq z) - \Phi(z)| \\ &\quad \leq \frac{9c_0^2 \mathrm{E}X_1^2}{(1+|z|)^2 n\sigma^2} + \frac{C}{(1+|z|)^p}\left(\frac{c_0\|X_1\|_2}{\sqrt{n}\sigma} + \frac{\mathrm{E}|g(X_1)|^p}{n^{(p-2)/2}\sigma^p}\right). \end{aligned} \tag{3.11}$$

## 4. An example

In this section we give an example to show that the bound of (2.4) in Theorem 2.1 is achievable. Moreover, the term $\sum \mathrm{E}|g_i(X_i)(\Delta - \Delta_i)|$ in (2.4) cannot be dropped. The example also provides a counterexample to a result of Shorack [22] and of Bolthausen and Götze [5]. We refer to Chen and Shao [10] for a detailed proof.

***Example 4.1.*** Let $X_1, \ldots, X_n$ be independent normally distributed random variables with mean zero and variance $1/n$. Define

$$W = \sum_{i=1}^n X_i, \qquad T := T_\varepsilon = W - \varepsilon|W|^{-1/2} + \varepsilon c_0 \quad \text{and} \quad \Delta = T - W = -\varepsilon|W|^{-1/2} + \varepsilon c_0,$$

where $c_0 = \mathrm{E}(|W|^{-1/2}) = \sqrt{2/\pi}\int_0^\infty x^{-1/2}\mathrm{e}^{-x^2/2}\,\mathrm{d}x$. Let $\{\hat{X}_i, 1 \leq i \leq n\}$ be an independent copy of $\{X_i, 1 \leq i \leq n\}$ and define

$$\alpha = \frac{1}{n}\sum_{i=1}^n \mathrm{E}|\Delta(X_1, \ldots, X_i, \ldots, X_n) - \Delta(X_1, \ldots, \hat{X}_i, \ldots, X_n)|. \tag{4.1}$$

Then $\mathrm{E}T = 0$ and for $0 < \varepsilon < 1/64$ and $n \geq (1/\varepsilon)^4$,

$$P(T \leq \varepsilon\, c_0) - \Phi(\varepsilon c_0) \geq \varepsilon^{2/3}/6, \tag{4.2}$$

$$\mathrm{E}|W\Delta| + \mathrm{E}|\Delta| \leq 7\varepsilon, \tag{4.3}$$

$$\mathrm{E}|\Delta| + \sum_{i=1}^n \mathrm{E}|X_i|^3 + \sqrt{\alpha} \leq C\varepsilon, \tag{4.4}$$

where $C$ is an absolute constant.

Clearly, (4.2) implies that

$$\sup_z |P(T_\varepsilon \le z) - \Phi(z)| \ge \frac{\varepsilon^{2/3}}{6}. \tag{4.5}$$

A result from Shorack [22] (see Lemma 11.1.3, page 261) states that for any random variables $W$ and $\Delta$,

$$\sup_z |P(W + \Delta \le z) - \Phi(z)| \le \sup_z |P(W \le z) - \Phi(z)| + 4\mathrm{E}|W\Delta| + 4\mathrm{E}|\Delta|. \tag{4.6}$$

Another result, which is in Theorem 2 of Bolthausen and Götze [5], states that if $\mathrm{E}T = 0$, then

$$\sup_z |P(T \le z) - \Phi(z)| \le C\left(\mathrm{E}|\Delta| + \sum_{i=1}^n \mathrm{E}|g_i(X_i)|^3 + \sqrt{\alpha}\right), \tag{4.7}$$

where $C$ is an absolute constant and $\alpha$ is defined in (4.1).

In view of (4.3), (4.4) and (4.5), the results of Shorack and of Bolthausen and Götze can be shown to lead to a contradiction.

## 5. Proofs of main theorems

**Proof of Theorem 2.1.** Equation (2.5) follows from (2.4) and (1.4). When $\beta > 1/2$, (2.4) is trivial. For $\beta \le 1/2$, (2.4) is a consequence of (2.3) and Remark 2.2. Thus, we need only to prove (2.3). Note that

$$-P(z - |\Delta| \le W \le z) \le P(T \le z) - P(W \le z) \le P(z \le W \le z + |\Delta|). \tag{5.1}$$

It suffices to show that

$$P(z \le W \le z + |\Delta|) \le 4\delta + \mathrm{E}|W\Delta| + \sum_{i=1}^n \mathrm{E}|g_i(X_i)(\Delta - \Delta_i)| \tag{5.2}$$

and

$$P(z - |\Delta| \le W \le z) \le 4\delta + \mathrm{E}|W\Delta| + \sum_{i=1}^n \mathrm{E}|g_i(X_i)(\Delta - \Delta_i)|, \tag{5.3}$$

where $\delta$ satisfies (2.2). Let

$$f_\Delta(w) = \begin{cases} -|\Delta|/2 - \delta, & \text{for } w \le z - \delta, \\ w - \frac{1}{2}(2z + |\Delta|), & \text{for } z - \delta \le w \le z + |\Delta| + \delta, \\ |\Delta|/2 + \delta, & \text{for } w > z + |\Delta| + \delta. \end{cases} \tag{5.4}$$

Let

$$\xi_i = g_i(X_i), \qquad \hat{M}_i(t) = \xi_i\{I(-\xi_i \le t \le 0) - I(0 < t \le -\xi_i)\},$$
$$M_i(t) = \mathrm{E}\hat{M}_i(t), \qquad \hat{M}(t) = \sum_{i=1}^n \hat{M}_i(t), \qquad M(t) = \mathrm{E}\hat{M}(t).$$

Since $\xi_i$ and $f_{\Delta_i}(W - \xi_i)$ are independent for $1 \le i \le n$ and $E\xi_i = 0$, we have

$$\begin{aligned}
\mathrm{E}\{W f_\Delta(W)\} &= \sum_{1\le i\le n} \mathrm{E}\{\xi_i(f_\Delta(W) - f_\Delta(W-\xi_i))\} \\
&\quad + \sum_{1\le i\le n} \mathrm{E}\{\xi_i(f_\Delta(W-\xi_i) - f_{\Delta_i}(W-\xi_i))\} \\
&:= H_1 + H_2. \qquad (5.5)
\end{aligned}$$

Using the fact that $\hat{M}(t) \ge 0$ and $f'_\Delta(w) \ge 0$, we have

$$\begin{aligned}
H_1 &= \sum_{1\le i\le n} \mathrm{E}\Big\{\xi_i \int_{-\xi_i}^0 f'_\Delta(W+t)\,\mathrm{d}t\Big\} \\
&= \sum_{1\le i\le n} \mathrm{E}\Big\{\int_{-\infty}^{\infty} f'_\Delta(W+t)\hat{M}_i(t)\,\mathrm{d}t\Big\} \\
&= \mathrm{E}\Big\{\int_{-\infty}^{\infty} f'_\Delta(W+t)\hat{M}(t)\,\mathrm{d}t\Big\} \\
&\ge \mathrm{E}\Big\{\int_{|t|\le\delta} f'_\Delta(W+t)\hat{M}(t)\,\mathrm{d}t\Big\} \\
&\ge \mathrm{E}\Big\{I(z \le W \le z + |\Delta|)\int_{|t|\le\delta} \hat{M}(t)\,\mathrm{d}t\Big\} \\
&= \sum_{1\le i\le n} \mathrm{E}\{I(z \le W \le z+|\Delta|)|\xi_i|\min(\delta, |\xi_i|)\} \\
&\ge H_{1,1} - H_{1,2}, \qquad (5.6)
\end{aligned}$$

where

$$H_{1,1} = P(z \le W \le z + |\Delta|)\sum_{1\le i\le n} \mathrm{E}\eta_i,$$
$$H_{1,2} = \mathrm{E}\Big|\sum_{1\le i\le n} \eta_i - \mathrm{E}\eta_i\Big|, \qquad \eta_i = |\xi_i|\min(\delta, |\xi_i|).$$

By (2.2),

$$\sum_{1\le i\le n} \mathrm{E}\eta_i \ge 1/2.$$

Hence

$$H_{1,1} \ge (1/2)P(z \le W \le z + |\Delta|). \tag{5.7}$$

By the Cauchy–Schwarz inequality,

$$\begin{aligned} H_{1,2} &\le \left(\mathrm{E}\left(\sum_{1\le i\le n} \eta_i - \mathrm{E}\eta_i\right)^2\right)^{1/2} \\ &\le \left(\sum_{1\le i\le n} \mathrm{E}\eta_i^2\right)^{1/2} \le \delta. \end{aligned} \tag{5.8}$$

As to $H_2$, it is easy to see that

$$|f_\Delta(w) - f_{\Delta_i}(w)| \le ||\Delta| - |\Delta_i||/2 \le |\Delta - \Delta_i|/2.$$

Hence

$$|H_2| \le (1/2)\sum_{i=1}^n \mathrm{E}|\xi_i(\Delta - \Delta_i)|. \tag{5.9}$$

Combining (5.5), (5.7), (5.8) and (5.9) yields

$$\begin{aligned} P(z \le W \le z + |\Delta|) &\le 2\left\{\mathrm{E}|Wf_\Delta(W)| + \delta + (1/2)\sum_{i=1}^n \mathrm{E}|\xi_i(\Delta - \Delta_i)|\right\} \\ &\le \mathrm{E}|W\Delta| + 2\delta\mathrm{E}|W| + 2\delta + \sum_{i=1}^n \mathrm{E}|\xi_i(\Delta - \Delta_i)| \\ &\le 4\delta + \mathrm{E}|W\Delta| + \sum_{i=1}^n \mathrm{E}|\xi_i(\Delta - \Delta_i)|. \end{aligned}$$

This proves (5.2). Similarly, one can prove (5.3) and hence Theorem 2.1. □

**Proof of Theorem 2.2.** First, we prove (2.9). For $|z| \le 4$, (2.9) holds by (2.5). For $|z| > 4$, consider two cases.

*Case 1*: $\sum_{i=1}^n \mathrm{E}|g_i(X_i)|^p > 1/2$. By the Rosenthal [20] inequality, we have

$$\begin{aligned} P(|W| > (|z| - 2)/3) &\le P(|W| > |z|/6) \le (|z|/6)^{-p}\mathrm{E}|W|^p \\ &\le C(|z| + 1)^{-p}\left\{\left(\sum_{i=1}^n \mathrm{E}g_i^2(X_i)\right)^{p/2} + \sum_{i=1}^n \mathrm{E}|g_i(X_i)|^p\right\} \end{aligned}$$

$$\le C(|z|+1)^{-p}\sum_{i=1}^{n} \mathrm{E}|g_i(X_i)|^p. \tag{5.10}$$

Hence

$$\begin{aligned}|P(T\le z)-\Phi(z)| &\le P(|\Delta|>(|z|+1)/3)+P(|W|>(|z|-2)/3)+P(|N(0,1)|>|z|)\\ &\le P(|\Delta|>(|z|+1)/3)+C(|z|+1)^{-p}\sum_{i=1}^{n}\mathrm{E}|g_i(X_i)|^p,\end{aligned}$$

which shows that (2.9) holds.

*Case 2*: $\sum_{i=1}^{n}\mathrm{E}|g_i(X_i)|^p\le 1/2$. Similar to (5.10), we have

$$\begin{aligned}P(|W-g_i(X_i)|>(|z|-2)/3) &\le C(|z|+1)^{-p}\left\{\left(\sum_{j=1}^{n}\mathrm{E}g_j^2(X_j)\right)^{p/2}+\sum_{j=1}^{n}\mathrm{E}|g_j(X_j)|^p\right\}\\ &\le C(|z|+1)^{-p}\end{aligned}$$

and hence

$$\begin{aligned}\gamma_z &\le P(|\Delta|>(|z|+1)/3)+\sum_{i=1}^{n}((|z|+1)/3)^{-p}\mathrm{E}|g_i(X_i)|^p+\sum_{i=1}^{n}C(|z|+1)^{-p}\mathrm{E}|g_i(X_i)|^p\\ &\le P(|\Delta|>(|z|+1)/3)+C(|z|+1)^{-p}\sum_{i=1}^{n}\mathrm{E}|g_i(X_i)|^p.\end{aligned}$$

By Remark 2.1, we can choose

$$\begin{aligned}\delta &= \left(\frac{2(p-2)^{p-2}}{(p-1)^{p-1}}\sum_{i=1}^{n}\mathrm{E}|g_i(X_i)|^p\right)^{1/(p-2)}\\ &\le \frac{2(p-2)^{p-2}}{(p-1)^{p-1}}\sum_{i=1}^{n}\mathrm{E}|g_i(X_i)|^p.\end{aligned}$$

Combining the above inequalities with (2.6) and the non-uniform Berry–Esseen bound for independent random variables yields (2.9).

Next we prove (2.6). The main idea of the proof is first to truncate $g_i(X_i)$ and then adopt the proof of Theorem 2.1 to the truncated sum. Without loss of generality, assume $z\ge 0$, because we can simply apply the result to $-T$. By (5.1), it suffices to show that

$$P(z-|\Delta|\le W\le z)\le\gamma_z+\mathrm{e}^{-z/3}\tau \tag{5.11}$$

and

$$P(z\le W\le z+|\Delta|)\le\gamma_z+\mathrm{e}^{-z/3}\tau. \tag{5.12}$$

Since the proof of (5.12) is similar to that of (5.11), we prove only (5.11). It is easy to see that

$$P(z-|\Delta| \le W \le z) \le P(|\Delta| > (z+1)/3) + P(z-|\Delta| \le W \le z, |\Delta| \le (z+1)/3).$$

Now (5.11) follows directly by Lemmas 5.1 and 5.2 below. This completes the proof of Theorem 2.2. □

**Lemma 5.1.** *Let*

$$\xi_i = g_i(X_i), \qquad \bar{\xi}_i = \xi_i I(\xi_i \le 1), \qquad \bar{W} = \sum_{i=1}^n \bar{\xi}_i.$$

*Then*

$$\begin{aligned} &P(z-|\Delta| \le W \le z, |\Delta| \le (z+1)/3) \\ &\quad \le P(z-|\Delta| \le \bar{W} \le z, |\Delta| \le (z+1)/3) \\ &\qquad + \sum_{i=1}^n P(\xi_i > (z+1)/3) + \sum_{i=1}^n P(W-\xi_i > (z-2)/3)P(|\xi_i| > 1). \end{aligned} \tag{5.13}$$

**Proof.** We have

$$\begin{aligned} &P(z-|\Delta| \le W \le z, |\Delta| \le (z+1)/3) \\ &\le P\Big(z-|\Delta| \le W \le z, |\Delta| \le (z+1)/3, \max_{1\le i\le n} |\xi_i| \le 1\Big) \\ &\quad + P\Big(z-|\Delta| \le W \le z, |\Delta| \le (z+1)/3, \max_{1\le i\le n} |\xi_i| > 1\Big) \\ &\le P(z-|\Delta| \le \bar{W} \le z, |\Delta| \le (z+1)/3) + \sum_{i=1}^n P(W > (2z-1)/3, |\xi_i| > 1) \end{aligned}$$

and

$$\begin{aligned} &\sum_{i=1}^n P(W > (2z-1)/3, |\xi_i| > 1) \\ &\quad \le \sum_{i=1}^n P(\xi_i > (z+1)/3) + \sum_{i=1}^n P(W > (2z-1)/3, \xi_i \le (z+1)/3, |\xi_i| > 1) \\ &\quad \le \sum_{i=1}^n P(\xi_i > (z+1)/3) + \sum_{i=1}^n P(W - \xi_i > (z-2)/3, |\xi_i| > 1) \\ &\quad = \sum_{i=1}^n P(\xi_i > (z+1)/3) + \sum_{i=1}^n P(W - \xi_i > (z-2)/3)P(|\xi_i| > 1), \end{aligned}$$

as desired. □

**Lemma 5.2.** *We have*

$$P(z - |\Delta| \le \bar{W} \le z, |\Delta| \le (z+1)/3) \le \mathrm{e}^{-z/3}\tau. \tag{5.14}$$

**Proof.** Noting that $\mathrm{E}\bar{\xi}_i \le 0$, that $\mathrm{e}^s \le 1 + s + s^2(\mathrm{e}^a - 1 - a)a^{-2}$ for $s \le a$ and $a > 0$ and that $a\bar{\xi}_i \le a$, we have for $a > 0$,

$$\begin{aligned}
\mathrm{E}\mathrm{e}^{a\bar{W}} &= \prod_{i=1}^n \mathrm{E}\mathrm{e}^{a\bar{\xi}_i} \\
&\le \prod_{i=1}^n (1 + a\mathrm{E}\bar{\xi}_i + (\mathrm{e}^a - 1 - a)\mathrm{E}\bar{\xi}_i^2) \\
&\le \exp\left((\mathrm{e}^a - 1 - a)\sum_{i=1}^n \mathrm{E}\bar{\xi}_i^2\right) \\
&\le \exp\left((\mathrm{e}^a - 1 - a)\sum_{i=1}^n \mathrm{E}\xi_i^2\right) \\
&= \exp(\mathrm{e}^a - 1 - a).
\end{aligned} \tag{5.15}$$

In particular, we have $\mathrm{E}\mathrm{e}^{\bar{W}/2} \le \exp(\mathrm{e}^{1/2} - 1.5)$. If $\delta \ge 0.07$, then

$$\begin{aligned}
&P(z - |\Delta| \le \bar{W} \le z, |\Delta| \le (z+1)/3) \\
&\quad \le P(\bar{W} > (2z-1)/3) \le \mathrm{e}^{-z/3+1/6}\mathrm{E}\mathrm{e}^{\bar{W}/2} \\
&\quad \le \mathrm{e}^{-z/3}\exp(\mathrm{e}^{0.5} - 4/3) \le 1.38\mathrm{e}^{-z/3} \le 20\delta\mathrm{e}^{-z/3}.
\end{aligned}$$

This proves (5.14) when $\delta \ge 0.07$.

For $\delta < 0.07$, let

$$f_\Delta(w) = \begin{cases} 0, & \text{for } w \le z - |\Delta| - \delta, \\ \mathrm{e}^{w/2}(w - z + |\Delta| + \delta), & \text{for } z - |\Delta| - \delta \le w \le z + \delta, \\ \mathrm{e}^{w/2}(|\Delta| + 2\delta), & \text{for } w > z + \delta. \end{cases} \tag{5.16}$$

Put

$$\bar{M}_i(t) = \xi_i\{I(-\bar{\xi}_i \le t \le 0) - I(0 < t \le -\bar{\xi}_i)\}, \qquad \bar{M}(t) = \sum_{i=1}^n \bar{M}_i(t).$$

By (5.5) and similar to (5.6), we have

$$\mathrm{E}\{W f_\Delta(\bar{W})\} = \mathrm{E}\left\{\int_{-\infty}^{\infty} f'_\Delta(\bar{W} + t)\bar{M}(t)\,\mathrm{d}t\right\}$$

$$
\begin{aligned}
&\quad + \sum_{i=1}^n \mathrm{E}\{\xi_i(f_\Delta(\bar W - \bar\xi_i) - f_{\Delta_i}(\bar W - \bar\xi_i))\} \\
&:= G_1 + G_2. \qquad (5.17)
\end{aligned}
$$

It follows from the fact that $\bar M(t) \ge 0$, $f'_\Delta(w) \ge \mathrm{e}^{w/2}$ for $z - |\Delta| - \delta \le w \le z + \delta$ and $f'_\Delta(w) \ge 0$ for all $w$ that

$$
\begin{aligned}
G_1 &\ge \mathrm{E}\biggl\{\int_{|t|\le\delta} f'_\Delta(\bar W + t)\bar M(t)\,\mathrm{d}t\biggr\} \\
&\ge \mathrm{E}\biggl\{\mathrm{e}^{(\bar W-\delta)/2} I(z - |\Delta| \le \bar W \le z, |\Delta| \le (z+1)/3) \int_{|t|\le\delta} \bar M(t)\,\mathrm{d}t\biggr\} \\
&= \mathrm{E}\{\mathrm{e}^{(\bar W-\delta)/2} I(z - |\Delta| \le \bar W \le z, |\Delta| \le (z+1)/3)\} \int_{|t|\le\delta} \mathrm{E}\bar M(t)\,\mathrm{d}t \\
&\quad + \mathrm{E}\biggl\{\mathrm{e}^{(\bar W-\delta)/2} I(z - |\Delta| \le \bar W \le z, |\Delta| \le (z+1)/3) \int_{|t|\le\delta} (\bar M(t) - \mathrm{E}\bar M(t))\,\mathrm{d}t\biggr\} \\
&\ge G_{1,1} - G_{1,2}, \qquad (5.18)
\end{aligned}
$$

where

$$
\begin{aligned}
G_{1,1} &= \mathrm{e}^{z/3-1/6-0.035} P(z - |\Delta| \le \bar W \le z, |\Delta| \le (z+1)/3) \int_{|t|\le\delta} \mathrm{E}\bar M(t)\,\mathrm{d}t, \\
G_{1,2} &= \mathrm{E}\biggl\{\int_{|t|\le\delta} \mathrm{e}^{\bar W/2} |\bar M(t) - \mathrm{E}\bar M(t)|\,\mathrm{d}t\biggr\}.
\end{aligned}
$$

By (2.2) and the assumption that $\delta \le 0.07$,

$$
\begin{aligned}
\int_{|t|\le\delta} \mathrm{E}\bar M(t)\,\mathrm{d}t &= \sum_{i=1}^n \mathrm{E}|\xi_i| \min(\delta, |\bar\xi_i|) \\
&= \sum_{i=1}^n \mathrm{E}|\xi_i| \min(\delta, |\xi_i|) \ge 1/2.
\end{aligned}
$$

Hence

$$
G_{1,1} \ge (1/2)\mathrm{e}^{z/3-1/6-0.035} P(z - |\Delta| \le \bar W \le z, |\Delta| \le (z+1)/3). \qquad (5.19)
$$

By (5.15), we have $\mathrm{E}\mathrm{e}^{\bar W} \le \exp(\mathrm{e} - 2) < 2.06$. It follows from the Cauchy–Schwarz inequality that

$$
G_{1,2} \le 0.5 \int_{|t|\le\delta} (0.5\mathrm{E}\mathrm{e}^{\bar W} + 2\mathrm{E}|\bar M(t) - \mathrm{E}\bar M(t)|^2)\,\mathrm{d}t
$$

$$\leq 0.5\left\{2.06\delta + 2\sum_{i=1}^{n}\int_{|t|\leq\delta}\mathrm{E}\xi_i^2(I(-\bar{\xi}_i \leq t \leq 0) + I(0 < t \leq -\bar{\xi}_i))\,\mathrm{d}t\right\}$$
$$= 0.5\left\{2.06\delta + 2\sum_{i=1}^{n}\mathrm{E}\xi_i^2\min(\delta, |\bar{\xi}_i|)\right\}$$
$$\leq 0.5\left\{2.06\delta + 2\delta\sum_{i=1}^{n}\mathrm{E}\xi_i^2\right\} \leq 2.03\delta. \tag{5.20}$$

With regard to $G_2$, it is easy to see that

$$|f_\Delta(w) - f_{\Delta_i}(w)| \leq \mathrm{e}^{w/2}||\Delta| - |\Delta_i|| \leq \mathrm{e}^{w/2}|\Delta - \Delta_i|.$$

Hence, by the Hölder inequality, (5.15) and the assumption that $\xi_i$ and $\bar{W} - \bar{\xi}_i$ are independent,

$$|G_2| \leq \sum_{i=1}^{n}\mathrm{E}|\xi_i\mathrm{e}^{(\bar{W}-\bar{\xi}_i)/2}(\Delta - \Delta_i)|$$
$$\leq \sum_{i=1}^{n}(\mathrm{E}\xi_i^2\mathrm{e}^{\bar{W}-\bar{\xi}_i})^{1/2}(\mathrm{E}(\Delta - \Delta_i)^2)^{1/2}$$
$$= \sum_{i=1}^{n}(\mathrm{E}\xi_i^2\mathrm{E}\mathrm{e}^{\bar{W}-\bar{\xi}_i})^{1/2}\|\Delta - \Delta_i\|_2$$
$$\leq 1.44\sum_{i=1}^{n}\|\xi_i\|_2\|\Delta - \Delta_i\|_2. \tag{5.21}$$

Following the proof of (5.15) and by using $|\mathrm{e}^s - 1| \leq |s|(\mathrm{e}^a - 1)/a$ for $s \leq a$ and $a > 0$, we have

$$\mathrm{E}W^2\mathrm{e}^{\bar{W}} = \sum_{i=1}^{n}\mathrm{E}\xi_i^2\mathrm{e}^{\bar{\xi}_i}\mathrm{E}\mathrm{e}^{\bar{W}-\bar{\xi}_i} + \sum_{1\leq i\neq j\leq n}\mathrm{E}\xi_i(\mathrm{e}^{\bar{\xi}_i} - 1)\mathrm{E}\xi_j(\mathrm{e}^{\bar{\xi}_j} - 1)\mathrm{E}\mathrm{e}^{\bar{W}-\bar{\xi}_i-\bar{\xi}_j}$$
$$\leq 2.06e\sum_{i=1}^{n}\mathrm{E}\xi_i^2 + 2.06(\mathrm{e} - 1)^2\sum_{1\leq i\neq j\leq n}\mathrm{E}\xi_i^2\mathrm{E}\xi_j^2$$
$$\leq 2.06\mathrm{e} + 2.06(\mathrm{e} - 1)^2 < 3.42^2.$$

Thus, we obtain

$$\mathrm{E}\{Wf_\Delta(\bar{W})\} \leq \mathrm{E}|W|\mathrm{e}^{\bar{W}/2}(|\Delta| + 2\delta)$$
$$\leq \{\|\Delta\|_2 + 2\delta\}(\mathrm{E}(W^2\mathrm{e}^{\bar{W}}))^{1/2}$$
$$\leq 3.42(\|\Delta\|_2 + 2\delta). \tag{5.22}$$

Combining (5.17), (5.19), (5.20), (5.21) and (5.22) yields

$$\begin{aligned} &P(z-|\Delta|\le \bar{W}\le z, |\Delta|\le (z+1)/3)\\ &\quad\le 2\mathrm{e}^{-z/3+1/6+0.035}\left\{3.42(\|\Delta\|_2+2\delta)+2.03\delta+1.44\sum_{i=1}^n \|\xi_i\|_2\,\|\Delta-\Delta_i\|_2\right\}\\ &\quad\le \mathrm{e}^{-z/3}\left\{22\delta+8.5\|\Delta\|_2+3.6\sum_{i=1}^n\|\xi_i\|_2\|\Delta-\Delta_i\|_2\right\}\\ &\quad= \mathrm{e}^{-z/3}\tau. \end{aligned}$$

This proves (5.14). □

## Acknowledgements

The authors are grateful to Xuming He for his contribution to the construction of the example in Section 4. The authors thank two referees, an associate editor and the editor for their valuable comments. The research of L.H.Y. Chen was partially supported by National University of Singapore grant R-146-000-013-112, and the research of Q.M. Shao was partially supported by NSF grant DMS-0103487, National University of Singapore grant R-146-000-038-101 and HKUST grant DAG 05/06.SC27 and RGC 602206.